\documentclass[final,epsfig,amsfont]{article}
\usepackage{latexsym}
\usepackage[dvips]{graphicx}

\newtheorem{lemma}{Lemma}[section]

\newcommand{\qed}{\rule{2mm}{3mm}}

\newcommand{\eeq}{\end{equation}}
\newcommand{\beq}{\begin{equation}} \newcommand{\nuq}[1]{\label{#1} \eeq}
\newcommand{\ba}{\begin{array}} \newcommand{\ea}{\end{array}}

\newtheorem{teo}{Theorem} \newcommand{\bth}{\begin{teo}} \newcommand{\enth}{\end{teo}}

\begin{document}
 
\title{Dynamical Systems and Numerical Analysis: the
Study of Measures generated by Uncountable~I.F.S.}
\author{G.Mantica\thanks{International Center for
Non-linear and Complex Systems, Universit\`a dell'Insubria, Via
Vallegio 11, Como, and CNISM, unit\`a di Como, I.N.F.N. sezione di
Milano, Italy E-mail: $<$giorgio@uninsubria.it$>$.}}
\date{}
\maketitle

 \begin{abstract}
Measures generated by Iterated Function Systems composed of uncountably many one--dimensional affine maps are studied. We present numerical techniques as well as rigorous results that establish whether these measures are absolutely or singular continuous.
 \end{abstract}
{\em keywords: Iterated Function Systems, Singular Measures, Fourier Transform,
Invariant Measures \\ AMS subject classification 47B38, 28A80, 60B10 }

\section{Introduction}
In this paper I want to describe an example of the fruitful interplay
between the theory of dynamical systems and numerical analysis: I want to
show how theoretical questions such as singularity (or continuity) with
respect to Lebesgue of a dynamical measure can be attacked from a
numerical point of view, and viceversa how particular features in the
numerical analysis of measures can be properly explained by concepts in
the theory of dynamical systems.

Let us consider the iteration of
maps $\phi_\lambda : X \rightarrow X$ from a compact metric space $X$ to
itself, labelled by the variable $\lambda$ which belongs to a measure
space $\Lambda$, on which the probability measure $\sigma$ is given. This
process gives rise to what is called an Iterated Function System, or
I.F.S., \cite{hut,barn} with invariant measure $\mu$, which can be defined
as follows. Consider the transfer operator $T$ on the space $C(X)$ of
continuous functions on $X$, via
   \beq
      (T h)(x) = \int d \sigma(\lambda)   (h \circ \phi_{\lambda}) (x).
   \nuq{inva1}
Then, let $T^*$ be the adjoint operator in the space of regular Borel
measures on $X$.  An invariant measure of the I.F.S. is the fixed point of
$T^*$,
    \beq
       T^*(\mu) = \mu.
   \nuq{inva2}
Suitable hypotheses can be formulated in order for $\mu$ to be unique,
given the choice of the set of maps $\phi_\lambda$ and of the distribution
$\sigma$ \cite{mendi}. We are interested in the nature of $\mu$: is it of
pure type? If so, is it pure point, singular, or absolutely continuous?
What characteristics of $\phi_\lambda$ and $\sigma$ have importance in
this regard? This problem belongs to a classical topic of research in
dynamical systems, that looks for {\em a.c.i.m.}, that is, absolutely
continuous invariant measures. Indeed, it is my opinion that {\em singular continuous} measures are equally, if not more
interesting, in many respects.

The plan of this paper is to study this problem in a class of maps of the
real line. We will derive algorithms to reduce the quest for an invariant measure to a fixed point problem. If in addition there exists an attractive fixed point in a suitable space of densities, absolute continuity of the measure will follow.
We will also describe the Fourier
transform of the measure $\mu$, its Mellin transform and its Sobolev
dimension, that will also lead to numerical and theoretical methods to determine absolute continuity or singularity of the invariant measure.

 \section{Infinite Affine Iterated Function Systems} \label{afifs}

Iterated Function Systems have been originally introduced and studied in
\cite{hut,barn}, although in some form their history goes much further
back in time, see \cite{zyg}. They have become a versatile mathematical tool
with applications to image compression \cite{barnima,lato1}, quantum dynamics \cite{etna} and much more \cite{diacone}. In its simplest form, an I.F.S. is a finite
collection of contractive maps of a space $X$ into itself: when iterated
randomly, these maps produce a stochastic process in $X$ with invariant
measure $\mu$. Interesting results and applications are found for affine
one dimensional maps of the kind
    \beq
     \phi_{\delta,\beta} (x) = \delta  (x - \beta) + \beta.  \label{ev427x}
 \eeq
Each of these maps has the fixed point $\beta$ and contractions ratio
$0<\delta<1$. When dealing with finite  collections of affine
one-dimensional maps, the problem of constructing their Jacobi matrix, and
therefore of computing integrals with respect to these measures, has been
solved in \cite{cap,hans,laurie}.

Rather than working with a finite set of maps associated with the pairs
$(\delta,\beta)$, as done usually, in this paper we adopt a generalization by which $\delta$ is fixed to a positive constant value,
strictly less than one, while $\beta$ is free to vary into a finite interval
of the real line, that, without loss of generality, we may understand to
be $[-1,1]$. This choice is called an infinite, homogeneous I.F.S.. It has
been originally investigated by Elton and Yan \cite{elton}.

Using the maps $\phi_{\delta,\beta}$ we can define a stochastic process in $X=[-1,1]$ via the
following rule: given an initial point $x \in [-1,1]$, choose a value of
$\beta$ at random in $[-1,1]$, according to a distribution $\sigma(\beta)$
(whose support may contain an infinity of points) and apply
$\phi_{\delta,\beta}$ to map $x$ into $\phi_{\delta,\beta}(x)$. Iterate
the procedure. A general theorem, due to Mendivil \cite{mendi}, guarantees
that there exists a unique invariant measure $\mu$ for this stochastic
process. In addition, this measure can be found, probability one, by the
Cesaro average of atomic measures at the points $x_j$ of a trajectory of
the process: $\frac{1}{n} \sum \delta_{x_j} \rightarrow \mu$. In a rather
pictorial view, we may describe the process above by saying that the point
$x$ is the location of a predator in a rather peculiar {\em chase} of the
prey, located at the point $\beta$: the predator moves towards the prey,
but as soon as its distance from $\beta$ is reduced by a factor $\delta$
the prey disappears, to reappear instantaneously at a new location, and
the process is repeated. Therefore, $\sigma$ is the distribution of the
position of the prey, that, along with the value of $\delta$ yields the
distribution $\mu$ of the predator's positions.
The measure $\mu$ can be equivalently defined by eq. (\ref{inva2}):
its approximation properties and the related inverse
problems have been discussed in \cite{steve,gio1,forte,gio-nal}.

The transfer operator $T$ for infinite, affine homogeneous I.F.S. takes the
following form:
 \beq
    (T f)(x)  =      \int   d\sigma(\beta) \;
      f \circ \phi_{\delta,\beta} (x) =
      \int   d\sigma(\beta) f(\delta x +   \bar{\delta} \beta),
  \nuq{bala1}
where $f$ is any continuous function, and where, to simplify the notation,
we have introduced the symbol $\bar{\delta}:= 1 - \delta$, to be used throughout the paper. Eq.
(\ref{inva2}) now becomes:
  \beq
        \int d(T^* \mu)(x) f(x) := \int d\mu(x) (T f)(x)
            = \int d\mu(x) \int   d\sigma(\beta) f(\delta x +   \bar{\delta} \beta).
\nuq{bala2a}
This is equivalent to say that, if $\mu$ is invariant, then the equality
  \beq
     \int d\mu(x) f(x) = \int d\mu(x) (T f)(x)
    \nuq{bala2}
holds for any continuous function $f$. In this setting, our problem can
therefore be formulated by asking what characteristics of the constant
$\delta$ and of the measure $\sigma$ influence the nature of $\mu$. We
start attacking this question in the next section, in the case when
$\sigma$ is the Lebesgue measure.

\section{Test Case I: $\sigma$ is the Lebesgue Measure on [-1,1]}
\label{sec-lebes}
In this section, we consider a case that can be treated
analytically as well as numerically. It is defined by letting $\sigma$ be
the Lebesque measure on $[-1,1]$. By proving theoretically the right
answer for the spectral type of $\mu$, we shall be able to use it as a
test of the numerical techniques that we shall introduce in the sequel. We
first prove a general result
   \begin{lemma}   \label{lem-cont}
   Suppose that $\sigma$ is absolutely continuous with a bounded density. Then, so is $\mu$.
      \end{lemma}
{\em Proof}.  Apply the balance relation (\ref{bala2})  to $f(x) =
\chi_{B_\epsilon(s)}(x)$, the characteristic function of $B_\epsilon(s)$,
the ball of radius $\epsilon$ centered at $s$ (technically, eq.
(\ref{inva2}) holds also for summable functions), to get:
 \beq
     \mu(B_\epsilon(s)) = \int d\mu(x) \sigma(B_{\epsilon/\bar{\delta}}(\frac{s - \delta x} {\bar{\delta}})).
   \nuq{ball1}
If $\sigma$ has a bounded density, this means that the following limit
exists and is equal to $\rho(\sigma;y)$, the density of the measure
$\sigma$ at the point $y$:
  \beq
            \lim_{\epsilon \rightarrow 0} \frac{1}{2 \epsilon}   \sigma(B_\epsilon(y)) = \rho(\sigma;y) .
      \nuq{density1}
Using this definition, we can now  divide both sides of eq. (\ref{ball1})
by $2 \epsilon$ and take the limit for $\epsilon \rightarrow 0$. Thanks to
the dominated convergence theorem the limit can be taken inside the
integral at r.h.s. to obtain
    \beq
        \rho(\mu;s) = \frac{1}{\bar{\delta}} \int d\mu(x) \rho( \sigma;\frac{s - \delta x} {\bar{\delta}}).
      \nuq{dens0}
From this last relation the thesis easily follows. \qed

We now specialize the theory to the case when $\sigma$ is the Lebesgue
measure, that is to say, all possible values of $\beta$ in $[-1,1]$ are
equally probable.
 \bth Let $0<\delta<1$ and let $\sigma$ be given by
 $
  \rho(\sigma;\beta)  = \frac{1}{2} \chi_{[-1,1]}(\beta).    $ 
Then, the invariant measure $\mu$ is absolutely continuous, with a density
$\rho(\mu;x)$ that is infinitely differentiable and non-analytic.
 \label{gio1}
   \enth
{\em Proof}. The fact that $\mu$ has a bounded density follows from Lemma
\ref{lem-cont}. Furthermore, inserting the value of $
  \rho(\sigma;\beta)$ into eq. (\ref{dens0}), we obtain
   \beq
    \rho(\mu;x) = \frac{1}{2 \bar{\delta}}  \mu( [\frac{x - \bar{\delta}}{\delta},\frac{x + \bar{\delta}}{\delta}] ),  \nuq{dens2}
and, taking the derivative with respect to $x$,
    \beq
      \rho'(\mu;x) = \frac{1}{2 \delta \bar{\delta}} [ \rho(\mu;\frac{x - \bar{\delta}}{\delta}) -\rho(\mu; \frac{x + \bar{\delta}}{\delta})].
 \nuq{dens3x}
This equation can be iterated to show the existence of derivatives of
$\rho(\mu;x)$ of all orders. Furthermore, observe that, when taking
$x=-1$, both points $(\frac{-1 - \bar{\delta}}{\delta})$ and $(\frac{-1 +
\bar{\delta}}{\delta})$ are to the left of $-1$, that is, outside of the
support of $\mu$, so that all derivatives of $\rho(\mu;x)$ are null in
$x=-1$. Therefore, $\rho(\mu;x)$ cannot be analytic. \qed

In the next section, we will describe a numerical technique to compute the
density of this invariant measure.

 \section{A Density Mapping}
 \label{densma}
Let us suppose that the invariant measure $\mu$ of an infinite, affine,
homogeneous I.F.S. is absolutely continuous, with density $\rho(\mu;x)$.
We want to derive a numerical technique to compute this density. In the
case of Bernoulli I.F.S., to be discussed later in section \ref{sec-bern},
this idea has been proposed in \cite{borw}, without quantitative tests of
convergence. Observe for starters that the action of the adjoint operator
$T^*$ can be transferred on densities. Use $f(x) =
\chi_{B_\epsilon(s)}(x)$ in eq. (\ref{bala2a}) to get, for any probability
measure $\nu$:
 \beq
   (T^* \nu)(B_\epsilon(x)) = \int d\sigma(\beta) \nu(B_{\epsilon/{\delta}}(\frac{x - \bar{\delta} \beta}{{\delta}})).
     \label{mub2}  \eeq
Clearly, a result quite similar to Lemma \ref{lem-cont} holds:
 \begin{lemma}   \label{lem-cont2}
  Suppose that $\nu$ is absolutely continuous with a bounded density. Then, so is $T^* \nu$.
\end{lemma}
In fact, operating as in the proof of Lemma \ref{lem-cont} of the previous
section, one also finds the density of $T^* \nu$:
   \beq
    \rho(T^* \nu;x) = \frac{1}{{\delta}} \int d\sigma(\beta) \rho(\nu;\frac{x - \bar{\delta} \beta}{{\delta}}).
    \nuq{dens3}
Eq. (\ref{dens3}) is particularly suited for analytical computations. It
can also be implemented numerically, whenever a convenient representation
of $\rho$ is found, and integration with respect to $\sigma$ is
numerically feasible.
In this perspective, it is best to consider a Fourier space
representation. Take therefore $f(x)=e^{-i y x}$ in eq. (\ref{bala2a}) and
use the notation
    \beq
         \hat{\nu}(y) := \int d\nu(x) e^{-i y x}
     \nuq{ft1}
to indicate the Fourier transform of an arbitrary measure $\nu$, to get
  \beq
    \widehat{T^* \mu}(y) = \hat{\mu}(\delta y) \cdot     \hat{\sigma}(\bar{\delta} y).
      \nuq{ft2}
This is a rather crucial relation, already discussed in \cite{elton}, that
links the Fourier transforms of $\sigma$ and $\mu$. We shall make use of
it repeatedly. Observe also that when $\nu$ is absolutely continuous, that
is, $d\nu(x) = \rho(\nu;x) dx$, $\hat{\nu}(y)$ in eq. (\ref{ft1}) can
equally well be seen as the Fourier transform of the density
$\rho(\nu;x)$. Furthermore, since the
support of $\mu$ is enclosed in $[-1,1]$, we consider the Fourier
coefficients of $\rho(\mu;x)$ over the basis set  $e^{-i \pi k
x}$, with $k$ integer. Let us call them $c_k(\mu)$:
   \beq
   c_k(\mu) := \hat{\mu}(\pi k) :=     \int d \mu(x) e^{-i \pi k x} =
\int_{-1}^1 e^{-i \pi k x}    \rho(\mu;x) dx .
        \nuq{ft3}
Setting $y=\pi k$ yields the Fourier coefficients of $T^* \mu$ at l.h.s.
of eq. (\ref{ft2}).  Simple formal manipulations of the right hand side
prove the Lemma:
    \begin{lemma}   \label{lem-four}
 The Fourier coefficients of $T^* \mu$ depend linearly on those of $\mu$
via the following operator
        \beq
           c_k(T^* \mu) =  \hat{\sigma}(\bar{\delta} \pi k)    \sum_{j=-\infty}^\infty \mbox{sinc} \, (\pi (j - k \delta)) \;    c_j(\mu).
             \nuq{ft4}
    \end{lemma}
This lemma is theoretically inspiring, and, at the same time, it can be
turned into a computational procedure. From the theoretical side, consider
the following argument. Suppose that we start from an initial distribution
$\mu_0$ whose Fourier coefficients decay extremely fast for $|k|$ larger
than some value, say $k_0$. As it is well known, this means that $\mu_0$
is absolutely continuous, with a very regular density. We observe
from eq. (\ref{ft4}) that this distribution of coefficients may {\em
spread} under the action of $T^*$. In fact, two phenomena are competing,
in this regard: On the one hand, $c_k(T^* \mu_0)$ receives a contribution
from $c_j(\mu_0)$ that is maximum for $j \sim k \delta$, and therefore
tends to ``populate'' higher frequencies, that is, to decrease the
regularity of the densities. On the other hand, this contribution is
multiplied by $\hat{\sigma}(\bar{\delta} \pi k)$, which tends to zero as
$k$ grows, if $\sigma$ is sufficiently regular, and therefore its effect
is diminished. On this basis, a theoretical estimate can be carried out,
in line with the project described in the introduction: use numerical
analysis to inspire theoretical proofs. This estimate will be presented
elsewhere. Instead, we now proceed with numerical techniques and
experimentations, that illustrate this point.

We use eq. (\ref{ft4}) as in a fixed point method: that is, we define a
sequence of measures  $\mu_n := T^{*n} \mu_0$, with $\mu_0$ an arbitrary
starting measure and we iterate eq. (\ref{ft4}). This leads to
   \begin{itemize}
   \item[] {\bf Algorithm 1.}
  \item[0] Fix an integer size $M$ and a threshold $\eta$
  and compute the vector of   Fourier coefficients
  $\hat{\sigma}(\bar{\delta} \pi k)$, for $k   = -M,\ldots,M$.
    \item[1] Initialization: Put $n=0$. Choose a suitable density to define
     the initial measure $\mu_0$ and its Fourier coefficients $c_k(\mu_0)$,
     for for $k   = -M,\ldots,M$.
    \item[2] Iteration: given the Fourier coefficients of $\mu_n$,
     compute $c_k(T^* \mu_n) = c_k(\mu_{n+1})$ via eq. (\ref{ft4}),
      restricting the summation over $j$ to the range $[-M,M]$.
   \item[3] Control and termination. Compute a distance (chosen among various
   possibilities, see later) between $\mu_n$ and  $\mu_{n+1}$.
   If difference is less than the threshold $\eta$,   stop.
   Otherwise, augment $n$ to $n+1$ and loop back to 1.
\end{itemize}
In numerical experimentations we use as initial measure $\mu_0$ either
the uniform density of $[-1,1]$, given by $\rho(\mu_0;x) = \frac{1}{2}
\chi_{[-1,1]}(x)$, whose Fourier transform is readily computed, or the
Fourier transform of a Gaussian distribution centered at zero, of variance
$\sqrt{2}/\omega$, with $\omega<1$.

\begin{figure}
\centerline{\includegraphics[width=6cm,height=12cm,angle=270]{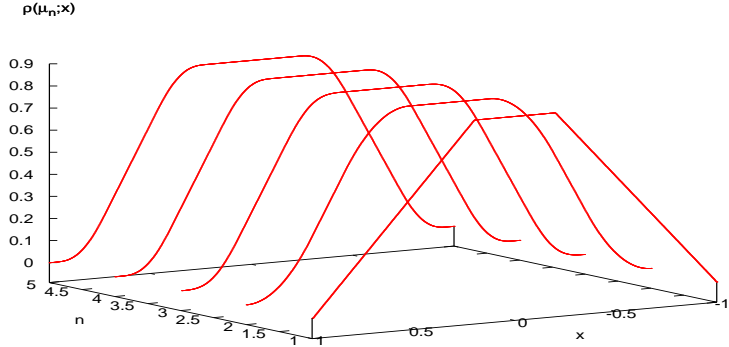}}
\caption{Densities $\rho(\mu_n;x)$ in the iteration of Algorithm 1, for
$n=1,\ldots,5$ when $\sigma$ is the Lebesgue measure. The initial density
$\rho(\mu_0;x)$ is constant.}
 \label{rofig1} \end{figure}

\begin{figure}
\centerline{\includegraphics[width=6cm,height=12cm,angle=270]{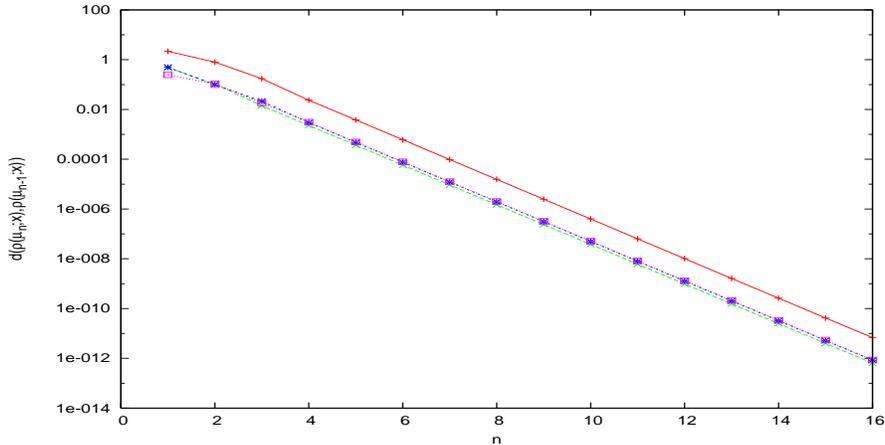}}
\caption{Distances $d_l(\rho(\mu_n;x),\rho(\mu_{n-1};x))$ between
successive densities in the iteration of Algorithm 1. Symbols for $l=4$
are red lines, crosses; $l=3$ blue lines, stars; $l=2$ green lines, stars;
and $l = 1$ magenta lines, squares.} \label{distafig2} \end{figure}

To test this technique we apply it first to the case where $\sigma$ is the
Lebesgue measure on $[-1,1]$, discussed in section \ref{sec-lebes}, for
which we have a precise theoretical result. In Fig. \ref{rofig1} we plot
the density of the measures obtained by the first five iterations of the
method, $\rho(\mu_n;x)$, $n=1,\ldots,5$, when and $\delta = \frac{2}{5}$
and $M=100$. The last two curves are practically undistinguishable,
indicating convergence. Indeed, a distance function is required to assess
this fact with rigour. Various choices are possible. In Fig.
\ref{distafig2} we plot all of the following distances, in the same case
of Fig. \ref{rofig1}, but on a larger number of iterations.
   \beq    d_l(\rho(\mu_n;x),\rho(\mu_{n-1};x)) =
     \left\{ \begin{array}{ll}
       \sum_{j=0}^M |c_j(\mu_n) - c_j(\mu_{n-1})|, & \; l = 1
  \\       \| \rho(\mu_n;x)- \rho(\mu_{n-1};x)\|_{1},& \; l = 2
    \\       \| \rho(\mu_n;x)- \rho(\mu_{n-1};x)\|_{\infty},& \; l = 3
      \\       \| \rho(\mu_n;x)- \rho(\mu_{n-1};x)\|_{BV},& \; l = 4
        \\
\end{array}    \right.
\nuq{dista1}
In all cases we observe exponential decrease of
$d_l(\rho(\mu_n;x),\rho(\mu_{n-1};x))$ with $n$, which implies convergence
to a density $\rho_M(\mu;x)$ that is an approximation of $\rho(\mu;x)$.
Validity of the approximation can be checked by increasing the value of
$M$, see the next Section and Sect. \ref{sec-sobo}, where a drastically
different case is encountered.

\section{Test Case II: $\sigma$ is a two-atoms Bernoulli Measure}
\label{sec-bern} Opposite to the previous choice of the measure $\sigma$, lies the Bernoulli measure, that is, an
atomic measure composed of two atoms. This is the first non--trivial case,
because one atom alone leads to the equality $\mu = \sigma$. Recall
the predator--prey interpretation: if the prey at $\beta$ does not move,
the predator $x$ converges to it, and the invariant distribution $\mu$ is
an atomic measure located at the same point. Let therefore $\sigma$ be
defined as the sum of two atomic measures located at $1$ and $-1$:
$\sigma(\beta) = \frac{1}{2}(\Delta_{-1}(\beta)+\Delta_1(\beta))$, or
 $
   \int f(\beta) d\sigma(\beta) = \frac{1}{2} (f(1) + f(-1)),
$
for any continuous $f$.

The prey appears with equal probability (this symmetry can be easily
changed, though) at either end of the interval. Contrary to what it might
seem, this problem is far from being trivial. It has a rich and classical
history, going back at least to Wiener, Wintner, Erd\"os  and others in
the 1930's. Present knowledge can be summed up in the
  \bth
 Let $\sigma(\beta) = \frac{1}{2}(\Delta_{-1}(\beta)+\Delta_1(\beta))$. For
any $\delta$, $\mu$ is of pure type. It is singular continuous, supported
on a Cantor set, if $\delta < \frac{1}{2}$. When $\delta = \frac{1}{2}$,
$\mu$ is the Lebesque measure on $[-1,1]$. There exist two constructive,
countable sets of values of $\delta> \frac{1}{2}$ for which $\mu$ is
singular continuous, and absolutely continuous, respectively. For almost
all $\delta> \frac{1}{2}$, $\mu$ is absolutely continuous. \enth

It is just the case to remark here that this theorem collect significant
results obtained along more than seventy years of research, see
\cite{borw,solom} for lists of references. It is easy to understand what
happens in the case $\delta < \frac{1}{2}$. Consider eq. (\ref{dens3}) and
let $\rho(\mu_0)$ be the uniform density on $[-1,1]$. It follows by direct
computation that $\rho(\mu_n)$ is a piece--wise constant function for any
$n$, which takes values $2^{-n-1} \delta^{-n}$ on a set of $2^n$ disjoint
intervals of equal length $2 \delta^n$, and zero otherwise. These
intervals constitute the usual generations in the hierarchical
construction of a Cantor set (take $\delta=\frac{1}{3}$
to obtain the classical, ternary Cantor set). Clearly, the sequence
$\rho(\mu_n)$ does not tend to any function and this is a remarkable example where measures converge while densities do not.

The bounded variation norm is particularly suited to illustrate this point: we have that
$\|\rho(\mu_n)\|_{BV} = \delta^{-n}$ and
$\|\rho(\mu_n)-\rho(\mu_{n-1})\|_{BV} = 2 \delta^{-n}(\frac{1}{2 \delta}
-1)$. Both quantities diverge as $n \rightarrow \infty$. This provides a
second test for algorithm I, that should {\em not} converge in this case. In numerical experiments we indeed observe
an exponential increase in the bounded variation norm
$\|\rho(\mu_n)\|_{BV}$, that only saturates because of the finite cardinality $2M+1$ of
the basis set. Quite obviously, increasing $M$ pushes the saturation
point to the right: higher and higher frequencies are required to describe
the densities $\rho(\mu_n)$. This much for Algorithm I: a better technique
will be described momentarily.

We need therefore to face numerically the possibility of $\mu$ not being
absolutely continuous, with slowly decaying (or not decaying at all!)
Fourier coefficients. In the next section we tackle this problem, first
theoretically and then numerically.

 \section{Fourier Transforms and Sobolev Dimension}
 \label{sec-sobo}
Let us therefore develop techniques for the case when $\mu$ may be
singular continuous. A few general results must be quoted at this point.
Consider the Mellin transform of a function $h$ defined on $[1,\infty]$:
   \begin{equation} \label{mel1}
   M_1 (h;z) := \int_1^{\infty}~y^{z-1} h(y) ~ dy.
   \end{equation}
The integral may diverge if $Re(z)$ is too large. The supremum of the set
of values of $Re(z)$ for which one has convergence is called the {\em
divergence abscissa} of the Mellin transform. Define the Sobolev dimension
of a measure $\mu$ as the divergence abscissa of the Mellin transform of
$h(y) = |\hat{\mu}(y)|^2$:
  \beq
d_s(\mu):= \sup \{  s \in {\bf R}  \; {  s.t. } \;
M_1(|\hat{\mu}(\cdot)|^2;s) < \infty \}.
    \nuq{sobo2}

Clearly, $d_s(\mu) \geq 0$. When $d_s(\mu)$ is less than, or equal to one,
it coincides with $D_2(\mu)$, the {\em correlation dimension} of the
measure $\mu$, a common quantity in the multifractal analysis of measures
\cite{pesin}. In this case, the leading asymptotic behavior for large $y$
of the Cesaro average of $h(y) = |\hat{\mu}(y)|^2$ is  $y^{-D_2}$,  see
\cite{stric,poin2} for a discussion of this and other asymptotic behaviors
of singular continuous measures. To the contrary, when $d_s(\mu) > 1$,
one knows that $\mu$ is absolutely continuous and its density $\rho(\mu)$
has fractional derivative of order $(d-1)/2$ in $L^2({\bf R})$ for all
$1<d< d_s(\mu)$. If in addition $d_s(\mu)>2$, then $\rho(\mu;x)$ is a
continuous function, the larger $d_s$ the higher its regularity. To sum
up, the regularity of a measure can be assessed by the study of the
asymptotic behavior of its Fourier transform.

To do this in our case, observe that equation (\ref{ft2}) can be iterated,
to get
   \beq
     \widehat{\mu_n}(y) := \widehat{T^{*n} \mu_0}(y) = \hat{\mu_0}(\delta^n y)
          \prod_{j=0}^{n-1}     \hat{\sigma}(\delta^j \bar{\delta} y).
     \nuq{ftr1}
Observe furthermore that, as $n$ tends to infinity, $\hat{\mu_0}(\delta^n
y)$ tends to $\hat{\mu_0}(0)=1$ for any $y$, so that this proves the
\begin{lemma}
 \label{lem-fprod}
 The Fourier transform of $\mu_n$ can be computed according to eq.
(\ref{ftr1}), and that of $\mu$ in the form of the infinite product
  \beq   \widehat{\mu}(y) = \prod_{j=0}^{\infty}     \hat{\sigma}(\delta^j \bar{\delta} y).
     \nuq{ftr3}
  \end{lemma}
Therefore, $\mu$ is an infinite convolution product of rescaled copies of
the measure $\sigma$, a fact well known when $\sigma$ is the Bernoulli
measure and $\widehat{\mu}(y)$ is an infinite product of trigonometric
functions. This observation is the basis of many theoretical
investigations. It also entails a numerical technique to compute
either $\widehat{\mu_n}(y)$ (just use eq. (\ref{ftr1})) or
$\widehat{\mu}(y)$:
  \begin{itemize}
   \item[] {\bf Algorithm 2.}
   \item[0] Fix a threshold $\eta>0$ and $y \in {\bf R}$.
     \item[1] Initialization: Put $n=0$ and let $\phi_0 = 1$.
        \item[2] Iteration: compute $\psi_n = \hat{\sigma}(\delta^n \bar{\delta}y)$.
        Update $\phi_{n+1} = \phi_{n} \psi_n$
        \item[3] Control and termination. Compute  $|\psi_n-1|$.
         If    this difference is less than the threshold $\eta$,    stop.
         Otherwise, augment $n$ to $n+1$ and loop back to 2.
         \end{itemize}

At this point, we can outline a two-fold strategy for detecting
numerically the continuity properties of a measure $\mu$ generated by
affine I.F.S.. In the first approach, we use Algorithm 2 to compute the
Fourier coefficients of $\mu$, and we look at their asymptotic behavior.

The second means is a refined density mapping: we fix a measure defined by
an initial density $\rho(\mu_0;x)$. Then, we compute the Fourier
coefficients $c_k(\mu_n)$ via eq. (\ref{ftr1}). Of course, while $\mu_n
\rightarrow \mu$ weakly because of Mendivil's theorem, there is no
guarantee that the densities $\rho(\mu_n;x)$ will tend to a function.
Therefore, we compute numerically the norms of $\rho(\mu_n;x)$
reconstructed from their Fourier coefficients, and the distances in eq.
(\ref{dista1}). We then look for either divergence of the bounded
variation norm (the $L^1$ norm must be conserved), or geometric
convergence of the distances. In the first case we conclude for
singularity of the measure $\mu$, in the second for absolute continuity.

We can test this approach by examining three cases of I.F.S. measures with
two maps, $\sigma = \frac{1}{2}(\Delta_{-1}+\Delta_1)$, and
$\delta>\frac{1}{2}$. The first two values of $\delta$ are chosen in the
two denumerable sets for which the spectral type is rigorously known:
absolutely continuous in the first case, a) $\delta = 2^{-1/2} \sim
0.7071067811865$ and singular continuous in the second case, b)
$\delta=1/p_1 \sim 0.7548776662467$, $p_1$ being a Pisot number. We also
consider the case c) $\delta=3/4=.75$ which is pretty close to b). The last
case is particularly interesting, since it is conjectured that for
rational values of $\delta$ larger than one half the measure $\mu$ is
absolutely continuous. In this specific case our analysis shows that the conjecture appears to be numerically validated.

\begin{figure}
\centerline{\includegraphics[width=6cm,height=12cm,angle=270]{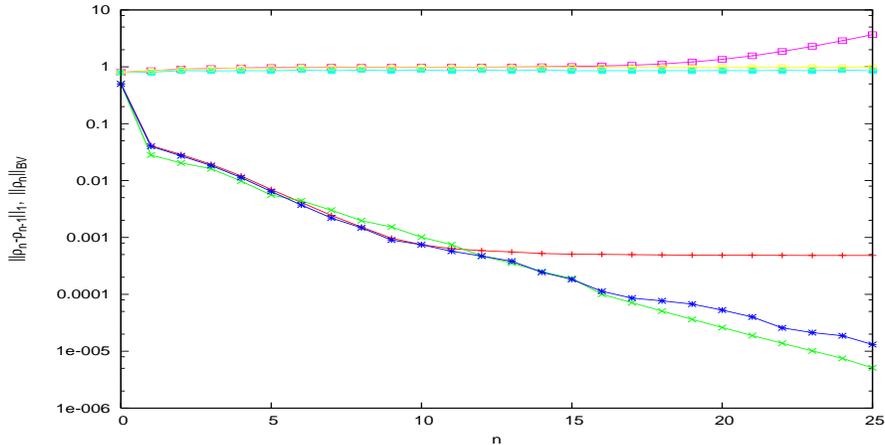}}
\caption{Bounded variation norms of densities $\rho_n$ for the three
I.F.S. described in the text: a) open squares, magenta; b) filled squares,
light blue; c) open circles, yellow. Also plotted are the $L^1$ distances
between successive densities: a) green crosses; b) red pluses; c) blue
asterisks.} \label{fig-mapd3} \end{figure}

Observe Fig. \ref{fig-mapd3}: the analysis of the $L^1$ distances
$\|\rho(\mu_n)-\rho(\mu_{n-1})\|_1$  versus $n$ shows exponential decay
(hence convergence in virtue of Cauchy criterion) for the known absolutely
continuous case b) and also for the conjectured case c). In case a) one
observes to the contrary (as expected) divergence of the bounded variation
norm of $\rho(\mu_n;x)$. Also notice that difference between cases b) and
c) can only be appreciated after about ten iterations.

\begin{figure}
\centerline{\includegraphics[width=6cm,height=12cm,angle=270]{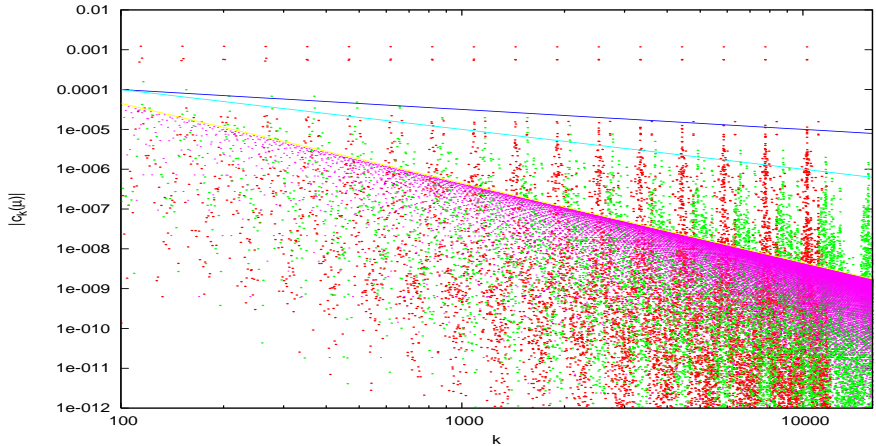}}
\caption{Fourier coefficients $|c_k(\mu)|$ vs. $k$ for the three I.F.S.
measures a) (red dots) b) (green dots) and c) (magenta dots) described in
Sect. \ref{sec-sobo}. Solid lines are drawn to mark the asymptotic
behaviors $k^{-1/2}$ (borderline for assessing absolute continuity),
$k^{-1}$ and $k^{-2}$.} \label{figfu1} \end{figure}

We also plot in fig. \ref{figfu1} the absolute values of the Fourier
coefficients $c_k(\mu)$ of the three invariant measures, versus $k$, in
double logarithmic scale. This confirms the the results obtained by the
analysis of the iterative algorithm. Coefficients for the case b) can also
be computed analytically.

\section{Test Case III: $\sigma$ is a Singular Continuous Measure} \label{sec-ifs}
We have seen so far two cases where $\sigma$ is a discrete measure or an absolutely continuous one. We want now to consider a
singular continuous measure $\sigma$. In order to be able to investigate this case
numerically, we need a fast and reliable way to compute the Fourier
transform of $\sigma$. But indeed we have already at hand a measure
$\sigma$ with these characteristics. Consider in fact the measure $\mu$
generated, for $\delta<\frac{1}{2}$, by the Bernoulli distribution,
section \ref{sec-bern} and \ref{sec-sobo}. It is singular continuous and
in addition its Fourier transform can be easily computed via Algorithm 2.
Therefore, by taking now this measure to be the {\em new} distribution
$\sigma$, we can apply Algorithm 2, in a recursive fashion.

The basic idea of this recursion is simple: suppose that the prey moves
with distribution $\sigma(\beta)=\mu^{(0)}(\beta)$, this latter being an
arbitrary measure. The predator is attracted to the prey according to the
affine maps with contraction ratio $\delta_1$: this determines
its invariant distribution $\mu^{(1)}$. In turn, this predator is hunted
by a second species: for this last animal,
$\sigma(\beta)=\mu^{(1)}(\beta)$ and $\delta_2$ is the contraction ratio:
its distribution is then $\mu^{(2)}$. Clearly, there is no limit to the
number of species in this food chain: this fact has an interesting
mathematical formulation that will be exploited elsewhere. For simplicity
we investigate in this section only a simple case: $\mu^{(0)}$ is the
Bernoulli measure, $\delta_1 < \frac{1}{2}$, and the choice
$\sigma(\beta)=\mu^{(1)}(\beta)$ leads to an I.F.S. with uncountably many
maps, whose fixed points populate a Cantor set.

\begin{figure}
\centerline{\includegraphics[width=6cm,height=12cm,angle=270]{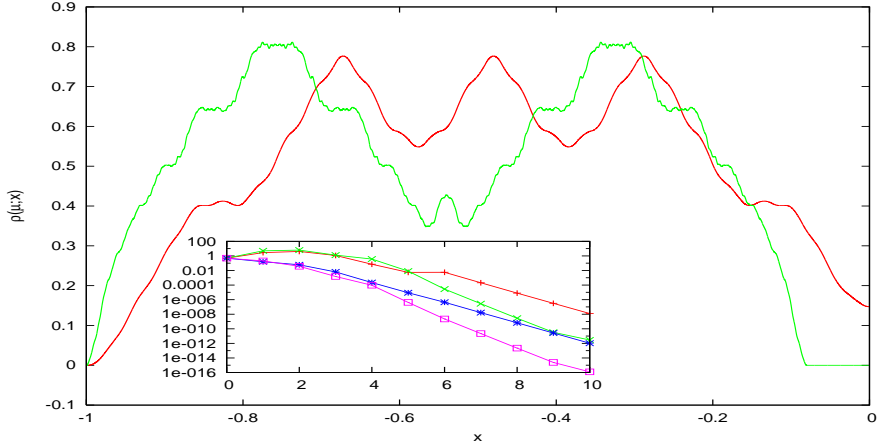}}
\caption{Densities $\rho(\mu^{(2)};x)$ of the two I.F.S. measures
described in the text. Labels are: $\delta_2=.1$ (green curve),
$\delta_2=.2$ (red curve). For symmetry, only half of the curves are
shown. In the inset, bounded variation and $L^1$ distances
$\|\rho(\mu_{n+1})-\rho(\mu_{n})\|$ versus iteration number $n$: labels
are $\delta_2=.1$, $BV$ (green curve, crosses), $L^1$ (magenta, squares);
$\delta_2=.2$, $BV$ (red, pluses), $L^1$ (blue, asterisks).
  }
 \label{figma3} \end{figure}

The theory developed so far enables us to examine numerically the
invariant measures of these I.F.S. As the typical case, we observe that
the measure $\mu^{(2)}$ is smoother than $\mu^{(1)}$. It can indeed be
absolutely continuous even when $\mu^{(1)}$ is singular continuous and
supported on a Cantor set. This case is pictured in Fig. \ref{figma3}, in
which $\mu^{(1)}$ is generated by an I.F.S. with
$\mu^{(0)}=\frac{1}{2}(\Delta_{-1}+\Delta_1)$ and $\delta_1=0.4$, and we
choose two values of the second contraction ratio: $\delta_2 = 0.1$ and
$\delta_2 = 0.2$. In both cases we observe absolute continuity of the
measure $\mu^{(2)}$, with a continuous density. The densities of the two
measures $\mu^{(2)}$ are plotted in Fig. \ref{figma3}. Observe that the
support of the measure $\mu^{(2)}$ includes that of $\mu^{(1)}$
\cite{gio1}. Indeed, a stricter ``inequality'' formula between them can be
derived. We also plot in the inset the geometric convergence of the $BV$
and $L^1$ distance functions (\ref{dista1}). The initial increase of the
bounded variation distance is due to the fact that we have employed here
an initial gaussian distribution (see sect. \ref{densma}).

Should one be
tempted to conjecture that $\mu^{(2)}$ is always absolutely continuous in
this setting, in the next and final section we shall meet a case when $\mu^{(2)}$ is singular continuous.

 \section{A computer assisted proof of absolute continuity} \label{cap}
In this final section, we show how a computer assisted proof of absolute
continuity of the measure $\mu$ can be developed in a specific case.

\begin{lemma}   \label{lem-cap}
Let the measure $\mu^{(0)}$ with support in $[-1,1]$ be given, and let
$\mu^{(1)}$ be the invariant measure generated by an affine, homogeneous
I.F.S. with contraction ratio $\delta$ and with distributions of fixed
points $\mu^{(0)}$. Let also $\mu^{(2)}$ be the invariant measure
generated by a second affine, homogeneous I.F.S. with contraction ratio
$\delta$ (the same as before) and with distributions of fixed points
$\mu^{(1)}$. Then, for any $N > 0$ and for any $t$, the inequality holds:
 \beq
   \log |\widehat{\mu^{(2)}} (\delta^{-N} t)| \leq N \{
  \log |\widehat{\mu^{(1)}} (\bar{\delta} t)| +
   \log |\widehat{\mu^{(0)}} (\delta^{-1} \bar{\delta}^2 t)| \}
  \nuq{eqcap1}
  \end{lemma}

{\em Proof.} Use eq. (\ref{ftr3}) to express $\log
\widehat{\mu^{(2)}}(y)$, and use the convenient notation $\Phi_j(y) :=
\log |\widehat{\mu^{(j)}}(y)|$. Gather terms (some algebra required) to
obtain
     \beq
  \Phi_2(\delta^{-N} t) = \sum_{l=0}^{N-1} (l+1) \Phi_0(\delta^{l-N}
  \bar{\delta}^2 t) + N \Phi_1(\bar{\delta} t) + \Phi_2(t).
     \nuq{ftr31}
Observe that all functions $\Phi_j$ are less than, or equal to, zero so
that by keeping only the $l=N-1$ term in the summation at r.h.s. of eq.
(\ref{ftr31}) and discarding also the term $\Phi_2(t)$ leads to the
inequality (\ref{eqcap1}). This is crucial for proving:
 \bth
 Let $\mu^{(j)}$, for $j=0,1,2$, be defined as in Lemma \ref{lem-cap}.
  Let also $\eta = - \sup \{\Phi_0(\delta^{-1} \bar{\delta}^2 t)
  + \Phi_1(\bar{\delta} t); \; t \in
[1,\delta^{-1}] \} $. Then, if $\eta > \frac{1}{2} \log(\delta^{-1})$ the
measure $\mu^{(2)}$ is absolutely continuous with an $L^2$ density, and if
$\eta > \log(\delta^{-1})$ the density $\rho(\mu;x)$ is continuous.
 \label{gio2}
   \enth

{\em Proof.} From Lemma \ref{lem-cap} it follows that for all $N>0$, and
for all $t \in [\delta^{-N},\delta^{-N-1}]$ the function $\Phi_2(t)$ is
less than $-\eta N$. One can therefore apply the Sobolev criterion
described in Sect. \ref{sec-sobo}. \qed

This theorem quickly turns into a computer--assisted proof of absolute
continuity. In fact, notice that the supremum required to compute $\eta$
is indeed the maximum of a continuous function defined on a finite
interval. As such, it is easily computable with sufficient precision to
establish whether the inequalities in Thm. \ref{gio2} are satisfied. For
the systems we have been considering in this paper this is easily
accomplished, by computing the Fourier transforms via algorithm 2. As just
an instance of this procedure, we can prove
 \bth
   Let $\mu^{(0)}=\frac{1}{2}(\Delta_{-1}+\Delta_1)$,
   and let $\mu^{(j)}$, $j=1,2$  defined as in Lemma \ref{lem-cap}.
 Then, if $\delta=\frac{2}{5}$ the measure $\mu^{(2)}$ is absolutely
continuous with a continuous density, and if $\delta=\frac{3}{10}$ it is
absolutely continuous with an
   $L^2$ density.
 \label{gio3}
   \enth
{\em Proof.} In the two cases quoted, one can easily compute that
$\eta/\log(\delta^{-1}) \simeq 1.5035$ and $\eta/\log(\delta^{-1}) \simeq
.96838$, respectively. \qed

It appears numerically that in addition to being a rigorous inequality,
formula (\ref{eqcap1}) is a rather precise estimate of the asymptotic
behavior of the Fourier transform of $\mu^{(2)}$. So, by letting now
$\delta=\frac{2}{10}$ we find that the value $\eta/\log(\delta^{-1})
\simeq 0.2218$ implies that the corresponding measure $\mu^{(2)}$ is
singular continuous. This conclusion is also confirmed by the iterative
analysis detailed in Sect. \ref{sec-sobo}.

\section{Conclusions}
\label{sec-conc} We have described numerical techniques to determine the
continuity type of measures generated by affine, homogeneous I.F.S. with
uncountably many maps. These techniques are inspired by the theory of
dynamical systems. In turn, they suggest methods to obtain rigorous, and
computer--assisted proofs: one should focus on the existence (or non--existence) of an attractive fixed point in a suitable space of densities, most conveniently under the bounded variation distance.  Alternatively, one must consider the decay rate of Fourier coefficients, also obtained via fixed point procedures.

We have sketched examples of both techniques, that have
helped us to clarify the fact that such I.F.S. measures can be of all
continuity types. According to the language adopted in this paper, we
have described the implications of the continuity of $\sigma$ on those of
$\mu$ in what we believe to be a rather detailed picture. We have also
shown applications of our methods to the classical problem of infinite
convolutions of Bernoulli measures, for which determining the continuity
type of a specific measure whose parameter does not belong to a known countable set of values is still an open problem.


\begin{thebibliography}{}



\bibitem{barn}
M. F. Barnsley and S. G. Demko, \newblock{Iterated function systems and
the global construction of fractals}. \newblock {\em Proc. R. Soc. London
A}, 399:243--275, 1985.

\bibitem{barnima}
M. F. Barnsley and L. P. Hurd,
 \newblock{Fractal Image Compression}.
  \newblock {\em A.K. Peters Ltd}, 1993.


\bibitem{steve}
D. Bessis and S. G. Demko, \newblock {Stable recovery of fractal measures
by polynomial sampling}. \newblock {\em Physica D} 47:427-438, 1991.

\bibitem{borw}
 J. Borwein, D. Bailey and R. Girgensohn,
 \newblock{Experimentation in Mathematics: Computational Paths to Discovery}.
 \newblock {\em A.K. Peters Ltd}, 2004.

\bibitem{diacone}
P. Diaconis and D. Freedman,
\newblock {Iterated random functions}.
\newblock {\em  SIAM Rev.},  41:45--76, 1999.

\bibitem{elton}
J. H. Elton and Z. Yan,
\newblock {Approximation of measures by Markov processes and
homogeneous affine iterated function systems}.
\newblock {\em  Constr. Appr.},   5:69--87, 1989.

\bibitem{hans}
H. J. Fischer,
\newblock {On generating orthogonal polynomials for self-similar measures}.
\newblock In N. Papamichael, St.
Ruscheweyh and E. B. Saff, editors, {\em Computational Methods and
Function Theory 1997}, 191--201, World Scientific Publishing Co. Pte. Ltd.
1999.

\bibitem{forte}
B. Forte and E.R. Vrscay,
\newblock {Solving the inverse problem for measures using iterated function systems: a new approach}.
\newblock {\em  Adv. Appl. Prob.},   27:800--820, 1995.


\bibitem{gio1}
C. R. Handy and G. Mantica,
\newblock{Inverse problems in fractal construction:
moment method solutions}.
\newblock{\em Physica D}, 43:17--36, 1990.

\bibitem{hut}
J. Hutchinson,
\newblock{Fractals and self-similarity}.
\newblock{\em Indiana J. Math.}, 30:713--747, 1981.


\bibitem{stric}
P. Janardhan, D. Rosenblum and R. S. Strichartz,
\newblock{Numerical experiments in Fourier asymptotics of Cantor measures and
wavelets}.
\newblock{\em Experiment. Math.}, 1:249--273, 1992.

\bibitem{lato1}
D. La Torre, E.R. Vrscay, M. Ebrahimi, M. F. Barnsley,
 \newblock{Measure-valued images, associated fractal transforms and the affine self-similarity of images},
 \newblock{\em SIAM Journal on Imaging Sciences}, 2:470--507, 2009.


\bibitem{laurie}
D. Laurie and J. de Villiers,
\newblock {Orthogonal polynomials for refinable linear functionals}.
\newblock {\em Math. Comp.},  75:1891--1903, 2006.

\bibitem{cap}
G. Mantica,
\newblock {A Stieltjes technique for computing
Jacobi matrices associated with singular measures}.
\newblock {\em Constr. Appr.},  12:509--530, 1996.



\bibitem{etna}
G. Mantica,
 \newblock {Fourier--Bessel functions of singular
continuous measures and their many asymptotics}. \newblock {\em E.T.N.A.} 25:409--430, 2006.

\bibitem{gio-nal}
G. Mantica,
 \newblock {Fractal measures and polynomial sampling: I.F.S.-Gaussian integration}.
 \newblock {\em Numer. Algor.} 45:269--281, 2007.

\bibitem{poin2}
G. Mantica and D. Guzzetti,
 \newblock {The asymptotic behaviour of the Fourier transform of
orthogonal polynomials II: Iterated Function Systems and Quantum
Mechanics}.
 \newblock {\em Ann. Henri Poincar\'e} 8:301--336, 2007.


\bibitem{mendi}
F. Mendivil,
 \newblock {A generalization of IFS with probabilities
to infinitely many maps}. \newblock {\em Rocky Mountain J. Math.}
28:1043--1051, 1998.


\bibitem{solom}
Y. Peres, K. Simon and B. Solomyak,
 \newblock {Absolute
continuity for random iterated function systems with overlaps}. \newblock
{\em J. London Math. Soc.} 74:739--756, 2006.


\bibitem{pesin}
Y. Pesin,
\newblock{Dimension Theory in Dynamical System:
Contemporary Views and Applications}.
\newblock {\em  Univ. Chicago Press}, 1996.


\bibitem{zyg}
A. Zygmund,
\newblock{Trigonometric series, Vols. I, II}.
\newblock {\em  Cambridge University Press}, 2002.

\end{thebibliography}
\end{document}